\documentclass[a4paper,12pt]{article}
\usepackage{amsmath}
\usepackage{amsfonts}
\usepackage{amssymb}
\usepackage{amsthm}
\usepackage{mathrsfs}
\usepackage[dvipsnames]{xcolor}

\newtheorem{theorem}{\bf Theorem}[section]
\newtheorem{lemma}[theorem]{\bf Lemma}

\newtheorem{definition}[theorem]{\bf Definition}
\newtheorem{remark}[theorem]{\bf Remark}
\newtheorem{proposition}[theorem]{ Proposition}

\newcommand{ \mis}{\textrm{meas}}  
\newcommand{\pr}[1]{
\begin{proof}
#1
\end{proof}
}
\newcommand{\es}[2]{
\begin{equation}\label{#1}
\begin{split}
#2
\end{split}
\end{equation}
}
\newcommand{\en}[1]{
\begin{equation*}
\begin{split}
#1
\end{split}
\end{equation*}
}

\usepackage[english]{babel}

\def \L {\mathscr{L}}
\def \LL {{\widetilde {\mathscr{L}}}}

\def \R {{\mathbb {R}}}
\def \N {{\mathbb {N}}}

\def \z {{\zeta}}

\def \tilde {\widetilde}

\def \H {\mathcal{Q}}
\begin{document}

\title{\bf Internal Schauder estimates \\ 
 for H\"ormander type equations \\ with Dini continuous source}

\author{
{\sc{Giovanna Citti\thanks{Dipartimento di Matematica, Universit\`{a} di Bologna, Piazza di Porta San Donato 5, 40126  Bologna, Italy. E-mail: giovanna.citti@unibo.it}}} 
\and 
{\sc {Bianca Stroffolini \thanks{ Dipartimento di Matematica, Universit\`{a} di Napoli Federico II, Via Cintia, 80126 Napoli, Italy.  E-mail: bstroffo@unina.it}}}}
\date{}
\maketitle
\begin{abstract}
\noindent We study the regularity properties of a general second order H\"ormander operator with Dini continous coefficients $a_{ij}$. Precisely if $X_0, X_1,\cdots X_m$  are smooth self adjoint  vector fields satisfying the  H\"ormander condition, we consider the linear operator in $\R^{N}$, with $N>m+1$:
\begin{equation*}
\L u := \sum_{i, j= 1}^{m} a_{ij} X_{i}X_{j}  u - X_0  u.
\end{equation*}
The vector field $X_0$ plays a role similar to the time derivative in a parabolic problem so that it is a vector of degree two.
We prove that, if  $f$ is a Dini continuous function, then the second order derivatives of the solution $u$ to the equation $\L u = f$ are Dini continuous functions as well. A key step in our proof is a Taylor formula in this anisotropic  setting, that we establish under minimal regularity assumptions.

\medskip

\noindent

\medskip

\end{abstract}

\setcounter{equation}{0}\setcounter{theorem}{0}
\section{Introduction}

In this article we consider a general second order operator  of the form 
\begin{equation}\label{e-Kolm-var}
\L u := \sum_{i, j= 1}^{m} a_{ij} X_iX_j   u - X_0  u,
\end{equation}
where the $a_{ij}$ is an uniformly  elliptic matrix, with Dini-continous coefficients, and $X_i$ are H\"ormander type vector fields and we study the local regularity of the  solution $u$ to $\L u = f$ when also $f$ is Dini-continuous.

When the matrix $a_{ij}$ is the identity, the operator reduces to a sum of squares of vector fields plus a drift term:
\begin{eqnarray} \label{operatorL}
   \LL :=\sum_{j=1}^m  X_j^2-X_0.
\end{eqnarray}
This class of operators has been introduced by Kolmogorov in \cite{Kolmogorov} and deeply studied by H\"{o}rmander's celebrated article \cite{Hormander}.  After that, the regularity theory has been widely developed in the seminal works by Folland \cite{Folland}, Folland and Stein \cite{FollandStein}, Rotschild and Stein \cite{RothschildStein}, Nagel, Stein and Wainger \cite{nagelstein}. 
In these papers the authors proved that the model operators for these problems are expressed in terms of vector fields invariant with respect to a Lie group structure. For this reason, a large literature have been developed under this extra assumption. 
We  refer to the recent monograph by Bonfiglioli, Lanconelli and Uguzzoni \cite{BLU} that contains an updated  description of this theory. If  $X_i=\partial_i$ for every $i=1, \cdots, m$, then the reduction to model Lie group has been studied by \cite{LanconelliPolidoro}, who were also able to provide an explicit expression of the fundamental solution, deeply simplifying the study of this special sub-class of operators. For this reason,  H\"older continuous Schauder estimates 
for this sub-class of operators have been investigated by Di Francesco-Polidoro in \cite{DiFrancescoPolidoro}, Lunardi in \cite{Lunardi}, Priola in \cite{Priola}, 
Wang-Zhang \cite{WangZhang}, and 
Biagi - Bramanti \cite{BiagiBramanti}, see also the recent papers \cite{LucertiniPascucci},\cite{BiagiBramantiStroffolini}.

In the elliptic and parabolic setting, regularity results are known also under the much weaker condition of Dini continuity, \cite{FabesSrokaWidman}, \cite{Kovats}, \cite{CaffarelliHuang}, 
\cite {Wang}. Precisely the condition 
 can be stated as follows: 

\begin{definition} \label{d-DC}
A function $f$  defined on an open set $H$ subset of a metric space  $(\R^{N+1}, d)$ is   Dini-continuous in $H$ if 
\begin{equation*}
\int_{0}^{1} \frac{\omega_f(r)}{r}dr < + \infty.
\end{equation*}
where $\omega_f$ is the modulus of continuity of $f$: 
\begin{equation} \label{eq-omegaf}
\omega_f(r) := \sup_{\substack{z, \zeta\in H \\ d(z,\zeta)<r}} |f(z)-f(\zeta)|.
\end{equation}
\end{definition}

The proof of \cite {Wang}, is much simpler, and it has been extended to prove internal Schauder estimates for hypoelliptic degenerate operators on the Heisenberg group by Wei, Jiang, and Wu  in \cite{WeiJiangWu}. 
In a different framework, Wang's method has been used by Bucur and Karakhanyan \cite{BucurKarakhanyan} in the study of fractional operators. 
Recently this approach has been applied by \cite{PRS} to a very special class of operators, where all the vector fields $X_j$ coincide with the partial derivatives $\partial_j,$ and the only non vanishing commutators are of the form $[\partial_j, X_0].$ In this case, the explicit expression of the fundamental solution is known, and a simple and direct procedure introduced in \cite{PRS} reduces the operator to a new one defined on an homogeneous Lie group, 
and the authors use these properties as a key ingredient of the proof of regularity of solutions, when the coefficients are only Dini continuous.

%
Our aim is to generalize these results to the large class of operators introduced in \cite{RothschildStein}, who proved very good local estimates of the fundamental solution even if exponential estimates are not known. Precisely, we only ask that $X_0, X_1,\cdots X_m$ do satisfy the H\"ormander condition, i.e. that the rank of the generated Lie algebra is maximum at every point. 
We will introduce a natural notion of degree, which assign degree 1 to the vector fields $ X_1,\cdots X_m$,   and degree two to the vector $X_0,$ since it plays the same role of the time derivative in the parabolic setting. In section 2 we also recall how to assign  to every tangent vector, a degree which  expresses the number of generators necessary to generate it from $X_0, X_1,\cdots X_m$. 
It is standard to associate a quasi distance $d_\L$ to the vector fields and their graduation.
This  notion of distance has  been first   studied by A. Nagel, E. M. Stein and S. Wainger \cite{nagelstein}, and it is expressed in terms of the exponential map.  Indeed there is a choice  $(X_j)_{j=0, \cdots N}$ of a basis of the tangent space at every point (see section 2) such that if  
\begin{equation}\label{espo}\zeta = \exp (\sum_j h_j X_j)(z), \end{equation}
then the distance between $z$ and $\zeta$
is expressed as 
\begin{equation}\label{dist}d_\L(z, \zeta) = \sum_j |h_j|^{1/deg(X_j)}, \end{equation}
where $deg$ denotes the degree just introduced, see \cite{RothschildStein}, \cite{Sanchez-Calle}, \cite{Morbidelli}. In addition to the standard exponential mapping, other maps have been introduced as composition of elementary translations along the vector fields (see \cite{nagelstein},  \cite{Danielli}). 
But only in \cite{Morbidelli}, the properties of these maps have been completely exploited. 
The exponential map introduced here  allows displacements only in the direction of the vector fields $X_0, \cdots, X_m$ and does not involve displacements along the direction of commutators, but it defines a metric equivalent to the classical one. 

Using the estimates of the fundamental solution, we  obtain  a first  estimate for the second order derivatives of the solutions, on the sphere of the metric. 
We define the cylinder $\H_R(z_0)$ using the geometric distance in this way:
$$\H_R(z_0)\colon= \{\z:d_\L(z_0, \z)<R\} .$$
\begin{proposition}\label{corollary}
Let $u$ be a solution to $\L u = 0$ in $\H_{R}(0)$, for $R \in ]0,1[$, then for any $X_i,X_j \in \lbrace X_1,...,X_m \rbrace$, there exists a constant $C$, only depending on $\lambda, \Lambda$ , such that
\begin{eqnarray*}
| X_i X_ j u |(z) \leq \frac{C}{R^2}\Vert u \Vert_{L^{\infty}(\H_{R}(0))},\quad z \in \H_{{\frac{R}{2}}}(0) \quad i, j= 1, \cdots m.
\end{eqnarray*}
\end{proposition}


We will also define class of $C^2_\L(\H_{1}(0))$ naturally associated to the vector fields and their graduation:

\begin{definition} \label{def-C2}
Let $\Omega$ be an open subset of $\R^{N+1}$. We say that a function $u$ belongs to $C^2_\L(\Omega)$ if $u$, its derivatives $X_i u, X_i X_j u$ ($i, j = 1, \dots, m$) and the Lie derivative $X_0 u$  are continuous functions in $\Omega$. 

We also require, for $i=1, \dots, m$, that
\begin{equation} \label{eq-Ymix}
 \lim_{s \to 0} \frac{X_i u (\text{\rm exp} (s X_0 )x, t-s) - X_i u(x,t)}{|s|^{1/2}} = 0,
\end{equation}
uniformly for every $(x,t) \in K$, where $K$ is a compact set $K \subset \Omega$. 
\end{definition}

We remark here that the condition \eqref{eq-Ymix} is minimal, since, if $X_0$ is the commutator of the other vector fields, then it is a necessary condition, as it has been proved in \cite{CittiManfredini}. Since $X_0$ plays the role of a second order derivative, condition \eqref{eq-Ymix} can be formally interpreted as a condition on the second order mixed derivative of the form $X_0^{1/2} X_i u$.

\begin{theorem} \label{th-2}
Let $\L$ be an operator in the form \eqref{e-Kolm-var}, and  $u$ be a classical solution to $\L u=f$. Assume  that $f$ and the coefficients $a_{ij}$, $i,j= 1, \dots, m$, are Dini continuous. Then for any points $z$ and $ \zeta \in \H_{\frac12}(0)$ the following holds:
\begin{equation*}
\begin{split}
&|X_iX_j u(z)-X_iX_j u(\zeta)| + |X_0 u(z)-X_0 u(\zeta)|\\
\le & c\Big( d \sup_{\H_{1}(0)} |u| + 
d \sup_{\H_{1}(0)}|f| + \int_0 ^d \frac{\omega_{f}(r)}{r}{dr}+ d \int_{d} ^1 \frac{\omega_{f}(r)}{r^2} dr\Big)  \\
+ & c \bigg( \sum_{i,j=1}^{m} \sup_{\H_{1}(0)} |\ \nabla^2_{x_i x_j} u| \bigg)
\Big(\int_0 ^d \frac{\omega_{a}(r)}{r}{dr}+ d \int_{d} ^1 \frac{\omega_{a}(r)}{r^2}{dr}\Big).
\end{split}
\end{equation*}
where $d=d_K (z, \zeta)$.
\end{theorem}

 One of the main ingredient of the proof is a new version of the anisotropic Taylor formula, for functions belonging to the natural  classes. Indeed we need to associated to a function $u$ around the point $z$ a polynomial function $P_z^2 u$ of order $2$ such that 
\begin{eqnarray}\label{TP}
u(\zeta)-P_z^2 u(\zeta)=o(d_\L(\z,z)^2) \quad \text{as} \quad \zeta \to z,
\end{eqnarray}
  Taylor development of this type associated to class of  vector fields is a well known result, which dates back to the  papers  \cite{RothschildStein}, and \cite{Folland}. In \cite{RothschildStein}, the function is of class $C^{\infty},$ and the anisotropic expansion was simply obtained by truncating the classical one, while in \cite{Folland} the condition of $C^{2,\alpha}$ was directly expressed by the existence of the approximating polynomial. When the regularity of the second order derivatives of $u$ is the very subject of the problem, the Taylor approximating result has to be established under minimal requirements on  $u$. For this reason, the requirements on the maximal derivative have been progressively lowered. Taylor formula for functions of class $C^2$ in the Euclidean sense has been proved by Arena, Caruso and Causa in \cite{ArenaCarusoCausa2010}.   
 Results which require that all the Euclidean second order derivatives of the function $u$ are H\"older continuous, are due to  Bonfiglioli \cite{Bonfiglioli} and to  Pagliarani, Pascucci and Pignotti \cite{PagliaPascPigno}. In \cite{PRS} a first result for functions of class $C^2_\L(\Omega)$ has been obtained 
 under the previously recalled very restrictive assumptions that the  vector fields $X_j$ coincide with $\partial_j$.  Here we prove the expansion in the full generality of the vector fields considered here, and only requiring that $u$ belongs to the space 
$C^2_\L(\Omega)$. It can be stated as following:

\begin{theorem} \label{taylor} 
Let $\L$ be the operator defined in \eqref{operatorL} , $\Omega$ be an open subset of $\R^{N+1}$ and let $u$ be a function in $C_{\L}^2(\Omega)$. For every $z \in \Omega$ we define the second order Taylor polynomial of $u$ around $z$ as
\begin{equation}\label{def-tay}
\begin{split}
 T^2_{z}u(\z) & := u(z) + \sum_{deg(X_i)=1} h_iX_i u(z) \\
& +\frac{1}{2} \sum_{i,j=1}^{m} h_i h_i X_i X_j  u(z)  + \sum_{deg(X_i)=2}h_i X_i u(z),
\end{split}
\end{equation}
for any $\z  \in \Omega$. Then condition \eqref{TP} is satisfied.  
\end{theorem}

\medskip

This paper is structured as follows. In Section 2, we recall some  basic notion of analysis in sub-Riemannian structures we need in our work. In particular, we recall the  properties about the fundamental solution of the constant coefficient operator $\LL$. In Section 3 we collect the  result regarding the  Taylor development of
a  function $u \in C^2_\L(\Omega)$.
In Section 4 we obtain the \emph{a priori} estimates of the second derivatives of the solutions stated in Proposition \ref{corollary}. Section 5 contains the proof of Theorem \ref{th-2} for the constant coefficient operator $\LL$, while Section 6 contains the proof for the general operator $\L$ with Dini-continuous coefficients. 

\setcounter{equation}{0}\setcounter{theorem}{0}
\section{The sub-Riemannian structure}

The vector fields $X_0, \cdots, X_m$ are H\"ormander type and self adjoint. They  generate a  Lie algebra of dimension $n=N+1$. We will define 
\begin{equation*}
\begin{split}
V_1 & = \text { span}(X_1, \cdots, X_m), 
 \quad V_2 = \text{ span } (V_1 \cup [V_1 , V_1] \cup \{X_0\}), \\ V_{j+1} & = \text{ span } (V_j \cup [V_1 , V_j] ) \quad  \text{ for every } \; j\geq 2.
\end{split}
\end{equation*}

Due to the H\"ormander condition, in a neighborhood of any point there exists 
a integer $s$ such that $V_s$ is the tangent space at that point. 
In particular, we assign $deg(X)=1$ for every $X\in  V_1$, $deg(x) = j$ if $X\in V_j$ and $X\notin V_{j-1}.$
The vector fields  will be expressed in coordinates as
$$X_i = \sum_{j=1}^n b_{ij}\partial_j,$$
where $b$ is of class $C^{\infty}.$
\subsection{Choice of the basis}\label{basis}
We will denote by $I=(I_0, \cdots I_N)$ any multi-index of length $N+1$. Next, we we will denote by $Y_1, \cdots Y_{\mathfrak {n}}$ an enumeration of all commutation of length at most $s$, where $s$ is defined as above as the smallest integer such that $V_s$ is the tangent plane  at the considered point. 
We select a basis of the space $Y_{I} $ as in \cite{Morbidelli} Theorem 2.1 (see also \cite{nagelstein}). 
Clearly, we can always  assume that the basis is ordered in such a way that, there exist $m_1, \cdots, m_s$, such that 
$Y_{I_1}, \cdots Y_{I_{m_1}}$ have degree 1, 
$Y_{I_{m_{i-1}+1}}=X_0, \cdots Y_{I_{m_{i}}}$ have degree $i$. 
Morbidelli in 
\cite{Morbidelli} introduced a modification of the exponential map defined only in terms of integral curves.

His definition is given by iteration. Let $S_1,\cdots, S_{\ell}$ a family of vector fields belonging to the family $X_0, X_1,\cdots, X_m$.
\begin{equation*}
\begin{split}
C_1(a; S_1)\colon &=  exp(a^{deg(S_1)} S_1) ;\\
C_2(a;S_1,S_2)\colon&= exp(-a^{deg(S_2)}S_2)exp(-a^{deg(S_1)}S_1)exp(a^{deg(S_2)}S_2)exp(a^{deg(S_1)}S_1)\\
\ldots\\
C_{\ell}(a; S_1,\cdots, S_{\ell})\colon&= C_{\ell-1}(a;S_2,\cdots, S_\ell)^{-1}exp(-a^{deg(S_1)}S_1)\\
&C_{\ell-1}
(a;S_2,\cdots, S_\ell) exp(a^{deg(S_1)}S_1)
\end{split}
\end{equation*}
Next, we denote by $d$ the sum of the degrees: $d=\sum_i deg(S_i)$.
For $\sigma>0$ small, 
\begin{equation*}
  exp^*(\sigma S_{[(1,\cdots, \ell]})  \colon=C_{\ell}(\sigma^{1/d};S_1,\cdots, S_\ell)
\end{equation*}

Then 
 $$E_I(z,h)= \exp^*(h_0 Y_0)\exp^*(h_1 Y_{I_1})...\exp^*(h_n Y_{I_n})(z) $$
 where $X_0=Y_0$. 
It has proved in \cite{Morbidelli} that $E_I(z,h)$ is a local diffeomorphism 
and $\frac{\partial E_I(z,h)}{\partial h_i} = Y_{I_i} + o(|h|),$
so that, by the choice of vector fields, the Jacobian of $E$ at the origin is
\begin{equation}\label{JacobianE}J E_I(z,0) = Id.\end{equation}
This expression induces a natural system of coordinates around the point $z$.   Indeed for every point $z$, and for every $\zeta= E_I(z,h),$ we will call $h= Log(\zeta)$ the coordinates of the point $\zeta$ around the point $z.$
Correspondingly, a distance $d$ can be defined  in terms of these coordinates: 
$$d_\L( z, \z)\colon=\sum_{i\in J} |h_i|^{\frac{1}{\deg(Y_i) }} $$
As proved in \cite{Morbidelli},  this distance is equivalent to the one defined in terms of the  standard exponential coordinates, recalled in \eqref{dist}.
We will call homogeneous dimension of the space in a neighborhood of the fixed point $$q = \sum_{i=0}^n deg(Y_{I_i}).$$
We define a family of dilation of a point $h= (h_0, \cdots h_n)$ on the Lie algebra: 
$$\delta_r(h)\colon=  (r^2 h_0, r h_1,  \cdots r^{deg(Y_{I_n})} h_n).$$
Accordingly, we have a family of local dilation  in a neighborhood of a fixed point  $z\in R^{N}$
\begin{equation} \label{dilations}
\delta_r(z) = E_I(0, \delta_r( Log(z))).
\end{equation}
As recalled in the introduction, one of the main result of \cite{RothschildStein}  and \cite{nagelstein} is that the operator $\LL$, defined in \eqref{e-Kolm-var}, admits  a fundamental solution  $\Gamma(x,y)$, which by definition satisfies:
$$u(x)=\int \Gamma(y,x)\LL(y) dy.$$
Moreover,
\begin{equation}\label{supgamma}
\begin{split}
|\Gamma(x,y)|\le c \frac{1}{d(x,y)^{q-2}}, \quad |X_j \Gamma(x,y)|\le c\frac{1}{d(x,y)^{q-2+deg(X_j)}},\\
 |X_i X_j \Gamma(x,y)|\le c\frac{1}{d(x,y)^{q-2+deg(X_j)+deg(X_i)}} 
 \end{split}
\end{equation}
\begin{remark}\label{ring}
We notice that, using the regularity of the vector fields, and, consequently, of the fundamental solution, we can obtain the following estimates on circular annuli:
\begin{equation*}
\sup_{\substack{z \in \H_{\frac{R}{2}}(0),\\  \zeta \in \H_{R}(0) \setminus \H_{\frac{3R}{4}}(0)}}\big\vert X_i \Gamma(z,\zeta)\big\vert \leq \frac{\tilde{C}}{R^{q-2+ deg(X_i)}}
\end{equation*}
for $\frac{3}{4}r\leq d_\L(x,y) \leq r.$ In addition, by composing a smooth  cut-off function with the distance, we can get:
\begin{equation*}
    \sup_{\H_{R}(0)  \setminus \H_{\frac{3R}{4}}(0)}\sum_{j= 1}^{m} |X_j  X_j \eta_{R}|\le \frac{C}{R^{\sum deg(X_j)}}
    \end{equation*}
\end{remark}

\setcounter{equation}{0}\setcounter{theorem}{0}
\section{Taylor development at order 2}
In this Section we prove Theorem \ref{taylor}. We generalize to our setting the procedure of Pagliarani, Pascucci and Pignotti \cite{PagliaPascPigno}, and \cite{PRS}. 
The main step of their work is to show that increments in the direction of the commutators can be obtained as displacements along suitable choice of increments along the horizontal vector fields. This fact has been established in full generality in Morbidelli in \cite{Morbidelli}, so that we will use his approach here. 
From now on, $Y_{I_j}$ will be a basis of the space, optimal in the
sense of section \ref{basis}.
\begin{proof}[Proof of Theorem \ref{taylor}]
 
By the results in \cite{Morbidelli}, we known that in a neighborhood of any fixed $z$, there exists an optimal basis $Y_I$ such that every $\z$ in a neighborhood of $z,$
there exist $h$ such that 
$$z = E_I(z,h).$$

Hence, it is sufficient to prove the results for 

 $$\z= \exp^*\Big(h_k Y_{I_k}\Big)(z),$$
with $h_k\to 0$ for every $k.$
 We will argue by finite induction with respect to the degree of  $Y_{I_k}$.
 
\medskip \noindent
\textbf{First step} Consider points 
$\z$ such that  
$$\z= \exp^*\Big(h_kY_{I_k}\Big)(z)$$
with $1\leq k \leq m,$ so that $deg(Y_{I_k})=1$, and 
$|h_{k}|\leq C$, for a fixed constant $C$.  Consequently, $\z= \exp\Big(h_kY_{I_k}\Big)(z),$ and  we can apply the Taylor formula in one variable to the curve $\gamma(h_k) = \exp\Big( h_kY_{I_k}\Big)(z),$ to obtain:
\begin{eqnarray}\label{opic}
u(\z) =u(z)+ h_kY_{I_k}u(z)  +\frac{h_k^2}{2} Y_{I_k}^2u(z)  + o(h_k^2),
\end{eqnarray}
which is the result in this case. 
\medskip

\medskip \noindent
\textbf{Inductive step}. 
We fix an integer $k$, assume that the thesis is true 
for every couple of points $z_1, \zeta_1$ such that 
$\zeta_1 = \exp^*(h_i Y_{I_i})(z_1) $, with deg$(Y_{I_i})<k$, and we prove the thesis for points
$$\zeta = \exp^*(h_tY_{I_t})(z),$$
with $deg(Y_{I_t}) =k$. Since  $deg(Y_{I_t})=k>1$ either $deg(Y_{I_t})=2$ and $Y_{I_t}= X_0,$
or there exist $I_i$ and $I_j$ such that $deg(i),$ $deg(j)< k,$ and 
$$Y_{I_t}= [Y_{I_i}, Y_{I_j}].$$

If $Y_{I_t}= X_0,$ we fix a  point $z, $ and consider points $\z$ in a neighborhood of $z$ such that 
 $$\z= \exp^*\Big(h_0Y_{I_0}\Big)(z)=  \exp\Big(h_0Y_{I_0}\Big)(z).$$
For $h_0$ in a neighborhood of $0$, this expression defines the general element of a curve $$\gamma: [0, C] \to R^n\quad \gamma(h_0) = \exp\Big(h_0 Y_{I_0}\Big)(z),$$
so that, we immediately obtain
\begin{eqnarray}\label{opic}
u(\z) =u(z)+ h_0 Y_{I_0}u(z) + o(h_0) = u(z)+ h_0 Y_{I_0}u(z) + o(d^2_\L(\zeta,z)).
\end{eqnarray}

In case $Y_{I_t} =[Y_{I_i}, Y_{I_j}],$
we assume by simplicity that $h_{t}>0$ and we define the points 
\en{
&z_1=\exp(  h_t^{deg(Y_{I_i})/ deg(Y_{I_t})  }  Y_{I_i} )(z),\\
&z_2=\exp(  h_t^{deg(Y_{I_j})/ deg(Y_{I_t})  }  Y_{I_j} )(z_1),\\
&z_3=\exp(- h_t^{deg(Y_{I_i})/ deg(Y_{I_t})  }  Y_{I_i} )(z_2)\\
&\zeta=\exp(- h_t^{deg(Y_{I_j})/ deg(Y_{I_t})  }  Y_{I_j} )(z_3).\\
}

If $deg(Y_{I_i})>2 $ and  $deg(Y_{I_j})>2$, then 
$$u(\zeta) = u(z) + o (h_t^{2/ deg(Y_{I_t})  } ) = u(z) + o(d^2_\L(z, \z) ).$$

If $deg(Y_{I_j})>2 $ and  $deg(Y_{I_i})\leq 2$, then 
 $$u(\zeta) = u(z_3) + (h_t^{2/ deg(Y_{I_t})  } ) = u(z_3) - u(z_2) + u(z_2) + o(h_t^{2/ deg(Y_{I_t})  } )= $$
$$= u(z_3) - u(z_2) + u(z_1) + o(h_t^{2/ deg(Y_{I_t})  }) =$$
$$= u(z_3) - u(z_2) - ( u(z)  -  u(z_1)) + u(z) + o(h_t^{2/ deg(Y_{I_t})  }) =$$
$$=T^2_{z_2}(z_3) - u(z_2) - (T^2_{z_1}(z) - u(z_1) ) + u(z) + o(h_t^{2/ deg(Y_{I_t})  }).$$

If  $deg(Y_{I_i})=1$
$$T^2_{z_2}(z_3) - u(z_2) =  h_t^{1/ deg(Y_{I_t})  } Y_{I_i}u(z_2) + \frac{h_t^{2/ deg(Y_{I_t})  }}{2} Y^2_{I_i}u(z_2)  $$
$$T^2_{z_1}(z) - u(z_1) = h_t^{1/ deg(Y_{I_t})  } Y_{I_i}u(z_1) + \frac{h_t^{2/ deg(Y_{I_t})  }}{2} Y^2_{I_i}u(z_1).$$
Then 
$$u(\zeta) = h_t^{1/ deg(Y_{I_t})  } (Y_{I_i}u(z_2) - Y_{I_i}u(z_1)) + \frac{h_t^{2/ deg(Y_{I_t})  }}{2} (Y^2_{I_i}u(z_2) - Y^2_{I_i}u(z_1))
= o(h_t^{2/ deg(Y_{I_t})  })$$

since $z_2 = z_1 + o(h_t^{2/ deg(Y_{I_t})  })$, the fact that $Y_{I_i}u$ is differentiable and $Y^2_{I_i}u$ is continuous. 

If  $deg(Y_{I_i})=2$, 
we have 
$$T^2_{z_2}(z_3) - u(z_2) =  h_t^{1/ deg(Y_{I_t})  } Y_{I_i}u(z_2)   $$
$$T^2_{z_1}(z) - u(z_1) =  h_t^{1/ deg(Y_{I_t})  } Y_{I_i}u(z_1) $$
Then 
$$u(\zeta) = h_t^{1/ deg(Y_{I_t})  } (Y_{I_i} u(z_2) - Y_{I_i} u(z_1)) 
= o(h_t^{2/ deg(Y_{I_t})  })$$
by the property on the fractional derivative that holds locally uniformly of the commutators proved in \cite{CittiManfredini} and the same property, required as an assumption in the direction $X_0$.

If $deg(Y_{I_j})\leq 2 $ and  $deg(Y_{I_i})\leq 2$ 
the proof is analogous.

This  concludes the proof of the inductive step.

\end{proof}

\section{$L^{\infty}$ estimates for  operators with constant coefficients}

In this Section we prove some a-priori estimates for the derivatives of a  solution to the  equation $\LL u = 0$  with constant coefficients and right-hand side equal to $0$. The solution is represented in terms of the fundamental solution $\Gamma$  of $\LL$, and the result is a direct consequence of the estimate of  $\Gamma$.
As a consequence, we obtain a simple mean-value formula for $u$.

\medskip

\medskip

In the sequel, we assume that all the eigenvalues of the constant matrix $A$ belong to some interval $[\lambda, \Lambda] \subset \R^+$. We are now in position to state our result. Recall that the complete enumeration of $X_i$ and their commutators will be denoted by $Y_j,$ with $j=1, \cdots \mathfrak{n}.$

\begin{proposition}\label{lem-apriori}
Let $u$ be a solution to $\LL u = 0$ in $\H_{R}(z_0)$, with $R \in ]0,1]$. Then 
\begin{eqnarray*}
| Y_j   u |(z) \leq \frac{C}{R^{deg(Y_j)}}\Vert u \Vert_{L^{\infty}(\H_{R}(z_0))},
\quad \text{for every} \quad z \in \H_{{\frac{R}{2}}}(z_0), \quad j=1,\ldots,\mathfrak{n},
\end{eqnarray*}
for some positive constant $C$ only depending on $\lambda, \Lambda$.
\end{proposition}
\begin{proof}
The proof parallels the one developed in \cite{PRS} with some little modifications.
The tools are in fact the a-priori estimates for the fundamental solutions and properties of the cut-off function.

Without loss of generality, we can assume $z_0=0$, since we will choose exponential coordinates in a neighorhood of an arbitrary point. Let $\eta_{R} \in C^{\infty}_{0}(\mathbb{R}^{N+1})$ be a cut-off function such that 
\begin{eqnarray}\label{cut-off}
\eta_{R}(z)=\chi(\Vert z\Vert_{\L}),
\end{eqnarray}
where $\chi \in C^{\infty}([0,+\infty),[0,1])$ is such that $\chi(s)=1$ if $s \leq \frac{3R}{4}$, $\chi(s)=0$ if $s \geq R$ and $|\chi'| \leq \frac{c}{R}$, $|\chi''|\leq \frac{c}{R^2}$. Then, for every $z \in \H_{R}(0)$ and for $i=1,\ldots,\mathfrak{n}$, there exists a constant $c$,  such that
\begin{eqnarray}\label{c_eta}
|Y_i \eta_{R}(z)| \leq \frac{c}{R^{deg(Y_i)}}.
\end{eqnarray} 
Consequently, for every $z \in \H_R(0)$ and $i,j=1,\ldots,m$, we have $|X_i X_j \eta_R (z)| \leq \frac{c}{R^2}$ and therefore we obtain a bound for the second order part of $\vert \LL \eta_R(z) \vert$. 

Since $\eta_R \equiv 1$ in $\in \H_{\frac{3R}{4}}(0)$, for every $z \in \H_{\frac{R}{2}}(0)$ we represent a solution $u$ to $\LL u=0$ as follows
\begin{align}\label{convolution}
u(z) = (\eta_{R}u)(z) =  \int_{\H_{R}(0)}(\Gamma(z,\cdot) \LL(\eta_{R}u))(\z)d\z.
\end{align}
Since $\LL=\sum_{i=1}^m X^2_i-X_0$ and $\LL u=0$ by assumption, \eqref{convolution} can be rewritten as
\begin{equation}\label{convolution2}
\begin{split}
u(z) = (\eta_{R}u)(z) &= + \int_{\H_{R}(0)}(\Gamma(z,\cdot) 
 (\sum_{j= 1}^{m} X_j  X_j \eta_{R})(\z) u(\z) d \z\\
&\quad- \int_{\H_{R}(0)}[\Gamma(z,\cdot) X_0 (\eta_{R})u](\z)d\z \\
&\quad+ 2\sum_{j= 1}^{m} \int_{\H_{R}(0)}[\Gamma(z,\cdot)(X_j u X_j\eta_{R}  ](\z)d\z.
\end{split}
\end{equation}
Integrating by parts the last integral in \eqref{convolution2}, we obtain, for every $z \in \H_{\frac{R}{2}}(0)$
\es{parts1}{
u(z) = (\eta_{R}u)(z) &=  + \int_{\H_{R}(0)}(\Gamma(z,\cdot) 
 (\sum_{j= 1}^{m} X_j  X_j \eta_{R})(\z) u(\z) d \z\\
&\quad- \int_{\H_{R}(0)}[\Gamma(z,\cdot) X_0 (\eta_{R})u](\z)d\z \\
&\quad+ 2 \sum_{j= 1}^{m}\int_{\H_{R}(0)}X_j(\Gamma(z,\cdot)X_j \eta_{R} )(\z) u(\z)d\z,
}
where we are differentiating w.r.t the variable $\zeta$. 

Since $\eta_R$ is constant in $\H_{\frac{3R}{4}}(0)$, all its derivatives vanish. After differentiating under the integral sign \eqref{parts1}, with respect to $Y_i$
for every $i=1,...,\mathfrak {n}$, we find
\en{
Y_i u(z) = Y_i(\eta_{R}u)(z) 
&= \int_{\H_{R}(0) \setminus \H_{\frac{3R}{4}}(0)} [Y_i\Gamma(z,\cdot) (\sum_{j= 1}^{m} X_j  X_j \eta_{R})(\z)  u](\z)d\z \\  
&\quad - \int_{\H_{R}(0) \setminus \H_{\frac{3R}{4}}(0)}[Y_i \Gamma(z,\cdot) X_0(\eta_{R}) u](\z)d\z \\ 
&\quad +2 \int_{\H_{R}(0) \setminus \H_{\frac{3R}{4}}(0)}\sum_{j= 1}^{m}  Y_i X_j(\Gamma(z,\cdot)X_j \eta_{R} )(\z) u(\z)d\z, \\
& =: \tilde{I}_{1}(z) + \tilde{I}_{2}(z)+\tilde{I}_{3}(z),
}

We estimate $\tilde{I}_{1}(z)$ and $ \tilde{I}_{2}(z)$, for $z \in \H_{\frac{R}{2}}(0)$. We have
\begin{align*}
\tilde{I}_{1}(z) &\leq \Vert u \Vert_{L^{\infty}(\H_{R}(0))}  \sup_{\H_{R}(0)  \setminus \H_{\frac{3R}{4}}(0)}\sum_{j= 1}^{m} |X_j  X_j \eta_{R}|\mis(\H_{R}(0)) \! \! \! \! \sup_{\substack{z \in \H_{\frac{R}{2}}(0),\\  \zeta \in \H_{R}(0) \setminus \H_{\frac{3R}{4}}(0)}}
\! \! \! \! \big\vert Y_i \Gamma(z,\zeta)\big\vert ,\\
\tilde{I}_{2}(z) &\leq \Vert u \Vert_{L^{\infty}(\H_{R}(0))}\sup_{\H_{R}(0) \setminus \H_{\frac{3R}{4}}(0)}\vert X_0(\eta_R)|\mis(\H_{R}(0)) \! \! \! \sup_{\substack{z \in \H_{\frac{R}{2}}(0),\\  \zeta \in \H_{R}(0) \setminus \H_{\frac{3R}{4}}(0)}} \! \! \! \big\vert Y_i \Gamma(z,\zeta)\big\vert.
\end{align*}
If $i= 0, \cdots \mathfrak{n}$, we  obtain
\begin{equation}\label{gamma_est}
\sup_{\substack{z \in \H_{\frac{R}{2}}(0),\\  \zeta \in \H_{R}(0) \setminus \H_{\frac{3R}{4}}(0)}}\big\vert Y_i \Gamma(z,\zeta)\big\vert \leq \frac{\tilde{C}}{R^{q-2+ deg(Y_i)}}.
\end{equation}
Moreover, by our choice of the cut-off function $\eta_R$, we have 
\begin{equation}\label{div_cut}
 |\sum_{j= 1}^{m} X_j  X_j(\eta_R)| \leq \frac{  c}{R^2} \quad \text{in} \quad \H_R(0).
\end{equation}
 Finally, combining inequalities \eqref{gamma_est} and \eqref{div_cut} with the fact that $\mis(\H_{R}(0))$ can be estimated by $C R^{q}$, for a suitable constant $C$ we obtain
\begin{eqnarray}\label{ine-above-2}
\tilde{I}_{1}(z) \leq \frac{C}{R^{deg(X_i)}}\Vert u \Vert_{L^{\infty}(\H_{R}(0))}, \quad z \in \H_{\frac{R}{2}}(0),
\end{eqnarray}
We now estimate $|X_0(\eta_{R})|$ in $\H_{R}(0) \setminus \H_{\frac{3R}{4}}(0)$. 
\begin{eqnarray}\label{add_term}
|X_0(\eta_{R})| \leq \frac{C'}{R^2}= \frac{C'}{R^{deg(X_0)}}, \quad \textrm{in $\H_{R}(0) \setminus \H_{\frac{3R}{4}}(0) $},
\end{eqnarray}
where $C'$ is a constant that only depends on the matrix $B$ and on the constant $c$ in \eqref{c_eta}.

Finally, using the fact that  $\mis(\H_{R}(0))$ can be estimated from below by $ C R^{q}$, together with \eqref{gamma_est} and \eqref{add_term}, we obtain 
\begin{eqnarray}\label{ine-above}
\tilde{I}_{2}(z) \leq \frac{C}{R^{deg(Y_i)}}\Vert u \Vert_{L^{\infty}(\H_{R}(0))}, \quad z \in \H_{\frac{R}{2}}(0),
\end{eqnarray}
where $C$ depends only on the constants $c$ and $\tilde{C}$ in \eqref{c_eta} and \eqref{gamma_est} .

By the same argument, we prove that, for a point $z \in \H_{\frac{R}{2}}(0) $, we have
\begin{eqnarray*}
\tilde{I}_{3}(z) \leq \Vert u \Vert_{L^{\infty}(\H_{R}(0))}\frac{c}{R}\ \mis(\H_{R}(0)) \! \! \! \! \! \! \! \! 
\sup_{\substack{z \in \H_{\frac{R}{2}}(0),\\  \zeta \in \H_{R}(0) \setminus \H_{\frac{3R}{4}}(0)}} \! \! \! \! \! \!
\big\vert \sum_{j= 1}^{m}  Y_i X_j(\Gamma(z,\cdot)\big\vert \\
\leq \frac{C}{R^{deg(Y_i)}}\Vert u \Vert_{L^{\infty}(\H_{R}(0))},
\end{eqnarray*}
where $C$ denotes once again a constant depending only on $c$, $\tilde{C}$ and $B$. Combining the inequality above with \eqref{ine-above-2} and \eqref{ine-above}, we finally obtain 
\begin{equation*}
\Vert Y_i u \Vert_{L^{\infty}(\H_{{\frac{R}{2}}}(0))} \leq \frac{C}{R^{deg(Y_i)}}\Vert u \Vert_{L^{\infty}(\H_{R}(0))}, \quad i=0,...,\mathfrak{n}.
\end{equation*}
\end{proof}

We state a result analogous to Proposition \ref{lem-apriori}, written in terms of the vector fields $X_1,\ldots,X_m,X_0$.
\begin{proposition}\label{corollary}
Let $u$ be a solution to $\LL u = 0$ in $\H_{R}(0)$, for $R \in ]0,1[$, then for any $X_i,X_j \in \lbrace X_1,...,X_m \rbrace$, there exists a constant $C$, only depending on $\lambda, \Lambda$ , such that
\begin{eqnarray*}
| X_i  u |(z) \leq \frac{C}{R^{deg(X_i)}}\Vert u \Vert_{L^{\infty}(\H_{R}(0))},\quad z \in \H_{{\frac{R}{2}}}(0),\quad i= 0, \cdots m \\
| X_i X_ j u |(z) \leq \frac{C}{R^2}\Vert u \Vert_{L^{\infty}(\H_{R}(0))},\quad z \in \H_{{\frac{R}{2}}}(0) \quad i, j= 1, \cdots m.
\end{eqnarray*}
\end{proposition}
\pr{The estimate of $X_1, \dots, X_m$ has been proved in Proposition \ref{lem-apriori}. The proof of the remaining estimates is obtained by reasoning as in Proposition \ref{lem-apriori}, and using estimates \eqref{supgamma} . We omit the details here.}

We now prove a mean value theorem for solutions $u$ to $\L u=0$ in cylinders $\H_R(\z)$.
\begin{proposition} [Scale invariant Lipschitz estimate] \label{mean-value-lem}
Let $\z$ be any point of $\R^{N+1}$, and let $u$ be a solution to $\L u=0$ in $\H_R(\z)$, with $R \in ]0,1]$. Then the following estimate holds
\begin{eqnarray}\label{lag}
\vert u (z)- u (\z) \vert \leq \frac{C}{R}d_\L(z,\z) \Vert u \Vert_{L^{\infty}(\H_{R}(\z))},
\end{eqnarray}
for every $z \in \H_\frac{R}{2}(\z)$. Here $C$ is a constant that only depends on $\lambda, \Lambda$ and on the matrix $B$.
\end{proposition}
\pr{
Since we will use exponential coordinates, it is not restrictive to assume $\z=0$, then we need to prove 
\begin{eqnarray*}
\vert u (z)- u (0) \vert \leq \frac{C}{R}\Vert z \Vert_\L \Vert u \Vert_{L^{\infty}(\H_{R}(0))}.
\end{eqnarray*}
Consider $z \in \H_{\frac{R}{2}}(0)$, Then there exists exponential coordinates $h_i$ such that 
$z = \exp(h_n Y_{I_n}) \cdots \exp( h_0 Y_{I_0})(0)$. By applying  the standard mean-value theorem we obtain
\begin{equation}\label{sum}
\begin{split}
|u(z)-u(0)| &= |u\Big(\exp(h_n Y_{I_n}) \cdots \exp( h_0 Y_{I_0})(0)\Big)-u(0)| =\\
&\Big| \sum_{i=1}^n u\Big(\exp(h_i Y_{I_i}) \cdots \exp( h_0 Y_{I_0})(0)\Big) - u\Big(\exp(h_{i-1} Y_{I_{i-1}}) \cdots \exp( h_0 Y_{I_0})(0)\Big)\Big|\\
&\leq \sum_{i=1}^{n} |h_i| \, |Y_{I_i} u\Big(\exp( \vartheta_i h_i Y_{I_i}) \cdots \exp( h_0 Y_{I_0})(0)\Big) |,
\end{split}
\end{equation}
where $\vartheta_i,\ldots,\vartheta_N,\vartheta \in ]0,1[$. For every $i=1,\ldots,n$, we have $$|h_i| \leq \Vert z \Vert_\L^{deg(Y_{I_i})} \leq R^{deg(Y_{I_i})},$$
and $\exp( \vartheta_i h_i Y_{I_i}) \cdots \exp( h_0 Y_{I_0})(0) \in \H_{\frac{R}{2}}(0)$. Then, by Proposition \ref{lem-apriori}, we find
\begin{eqnarray*}
|Y_{I_i} u \Big(\exp( \vartheta_i h_i Y_{I_i}) \cdots \exp( h_0 Y_{I_0})(0) \Big)| \leq \frac{c}{R^{deg(Y_{I_i})}} \Vert u \Vert_{L^{\infty}(\H_{R}(0))}.
\end{eqnarray*}
so that
\begin{eqnarray*}
|h_i| \, \Big|Y_{I_i} u \Big(\exp( \vartheta_i h_i Y_{I_i}) \cdots \exp( h_0 Y_{I_0})(0) \Big)\Big| \leq 
\frac{c}{R} \Vert z \Vert_\L \Vert u \Vert_{L^{\infty}(\H_{R}(0))}.
\end{eqnarray*}

The proof of the proposition can be obtained by combining the above estimates.
}

\begin{lemma}\label{lem-tec}
Let $\eta_R$ be the cut-off function introduced in \eqref{cut-off}. Then there exists a positive constant $C$, such that
\begin{eqnarray}\label{est-sec-gamma-const}
\bigg\vert  X_{I_i}X_{I_j} \int_{\H_R(z)}\Gamma(z,\z)\eta_R(\z)  d\z \bigg\vert \leq C,
\end{eqnarray}

for every $z \in \H_{\frac{R}{2}}(0)$, $R \in ]0,1]$ and for any $i,j=1,\ldots,m$.

\end{lemma}
\begin{proof}
The proof is a standard consequence of the definition of fundamental solution. We provide here the proof for reader convenience. (see also \cite{Stein}, or  \cite{DiFrancescoPolidoro}). We write the right-hand side of \eqref{est-sec-gamma-const} as 
\begin{equation}\label{2.111}
\begin{split}
X_{I_i}X_{I_j}&\int_{\H_R(z)} \Gamma(z, \z) \eta_R(\z) d\z \\&=X_{I_i}\int_{\H_R(z) }X_{I_j}\Gamma(z, \z) \eta_R(\z) d\z 
\end{split}
\end{equation}
We will need to write the derivative in terms of the variable $z$, denoted $X^z_{I_j}$ in terms of a derivative  with respect to $\z$, denoted $X^\z_{I_j}$. Indeed we can apply the following  formula, proved at the end of page 295 in \cite{RothschildStein}: 
$$X^z_{I_i} = \sum_{k} b_{i, k} X^\z_{I_k}$$
where $b_{j, k}$  is a function of degree greater or equal to $deg(X_{I_k}) - deg(X_{I_j})$. Consequently 

\begin{equation}\label{2.111}
\begin{split}
X_{I_i}X_{I_j}&\int_{\H_R(z)} \Gamma(z, \z) \eta_R(\z) d\z =\\
&=\int_{\H_R(z) }X_{I_i}X_{I_j}\Gamma(z, \z) \left(\eta_R(\z)-\eta_R(z)\right) d\z \\
&\quad - \eta_R(z) \int_{\H_R(z)} 
 \sum_k  X^\z_{I_k} b_{i k}  X_{I_j}\Gamma(z, \z) d\z \\
 &+ \eta_R(z) \int_{\partial \H_R(z)} 
 \sum_k   b_{i k}  X_{I_j}\Gamma(z, \z) \nu _{I_k}d\sigma (\z)=: I_1(z) + I_2(z) + I_3(z).
\end{split}
\end{equation}
By the definition of $\eta_R$, we have
\begin{eqnarray}\label{bound-eta}
0 \leq \eta_R \leq 1, \quad \eta_R(\z)-\eta_R(z) =0, \quad \forall \z \in \H_{\frac{3R}{4}}(z), z \in \H_{\frac{R}{2}}(0).
\end{eqnarray}
Thus, taking advantage of Remark \ref{ring}, we infer
\begin{equation}\label{est-const}
\begin{split}
&|I_1(z)| = \bigg\vert\int_{\H_R(z) }X_{I_i}X_{I_j}\Gamma(z, \z) \left(\eta_R(\z)-\eta_R(z)\right) d\z\bigg\vert\\
&\quad=\bigg\vert\int_{\H_R(z) \setminus \H_{\frac{3R}{4}(z)}}X_{I_i}X_{I_j}\Gamma(z, \z) \left(\eta_R(\z)-\eta_R(z)\right) d\z \bigg\vert \leq C
\end{split}
\end{equation}
Note that there exists a constant $C$ such that 
$$|I_2(z)| + |I_3(z)|\leq C,$$
since $X_{I_j}\Gamma$ is locally integrable and 
it has no singularities when $\z \in \partial \H_{R}(z)$.

\end{proof}

\setcounter{equation}{0}\setcounter{theorem}{0}
\section{H\"older estimates for  operators with constant coefficients}
We first prove a preliminary lemma, which is a straightforward consequence of the maximum principle.
\begin{lemma}\label{lemma-max}
Given $\varphi \in C(\partial \H_R(z_0))$ and $g \in C_b(\H_R(z_0)) $, we let $v$ be the solution to the following Dirichlet problem
\begin{equation*}
\left\{ \begin{array}{ll}
\L v=g,\quad &\textit{in $\H_R(z_0)$},\\
v=\varphi,\quad &\textit{in $\partial \H_R(z_0)$}.
\end{array} \right.
\end{equation*}
Then, the following holds
\begin{eqnarray}\label{max-principle}
\Vert v \Vert_{{L}^\infty(\H_R(z_0))}\leq \Vert \varphi \Vert_{{L}^\infty(\H_R(z_0))} + R^2 \Vert g \Vert_{{L}^\infty(\H_R(z_0))}.
\end{eqnarray}

\end{lemma}
\pr{
Assume that $z= E(z_0, h)$ 
We introduce the function $w(z): =2(h_0+ R^2)\Vert g \Vert_{{L}^\infty(\H_R(z_0))}+\Vert \varphi \Vert_{{L}^\infty(\H_R(z_0))}$ and we let $u:=v-w$.
Using the results in \cite{Morbidelli}, we see that $$X_0 h_0 = 1 + o(R).$$
Hence, for $R$ sufficiently small 
$$X_0 h_0 \in [\frac{1}{2}, \frac{3}{2}].$$

Clearly, $u$ satisfies $\L u = g+(2 + o(1))\Vert g \Vert_{{L}^\infty(\H_R(z_0))} \geq 0$ in $\H_R(z_0)$. Moreover, as $v \equiv \varphi$ on the boundary of $ \H_R(z_0)$, we have $u=\varphi-2(h_0 + R^2)\Vert g \Vert_{{L}^\infty(\H_R(z_0))}-\Vert \varphi \Vert_{{L}^\infty(\H_R(z_0))}\leq \varphi-\Vert \varphi \Vert_{{L}^\infty(\H_R(z_0))}\leq 0$ in $\partial \H_R(z_0)$. By the strong maximum principle, it follows that $u(x,t) \leq 0$ in $\H_R(z_0)$. 
Replacing $v$ by $-v$, estimate \eqref{max-principle} follows at once.
}

\pr{[Proof of Theorem \ref{th-2}]  We first consider the case of constant coefficients. We denote $\H_{k}=\H_{\varrho^{k}}(0)$, $\varrho=\frac{1}{2}$ and we consider the following sequence of Dirichlet problems:
\begin{equation}\label{dirichlet1}
\left\{ \begin{array}{ll}
\L u_{k}=f(0),\quad \textit{in $\H_{k}$}\\
u_{k}=u,\quad \textit{in $\partial \H_{k}$}
\end{array} \right.
\end{equation}
For any point $z=(x,t)$ satisfying $\Vert z \Vert_K \leq \frac{1}{2}$, we want to estimate the quantity
\begin{eqnarray*}
I(z) := |\partial^2 u(z) - \partial^2 u(0)|,
\end{eqnarray*}
where $\partial^2 u(z)$ stands for either $X_{I_i}X_{I_j} u(z)$, with $i,j=1,\ldots,m$, or $X_0 u(z)$.
To this end, we write $I$ as the sum of three terms:
\begin{align*}
I(z) &\leq |\partial^2 u_{k}(z) - \partial^2 u_{k}(0)| + |\partial^2 u_{k}(0) - \partial^2 u(0)| + \\
&\qquad\qquad\quad+|\partial^2 u(z) - \partial^2 u_{k}(z)| =: I_{1}(z)+I_{2}(z)+I_{3}(z).
\end{align*}
We first estimate $I_{2}$. Following \cite{Wang}, we prove that $\left(\partial^2 u_{k}(0)\right)_{k \in \N}$ is a Cauchy sequence and that its limit agrees with $\partial^2 u(0)$. The same assertion holds for $I_3$ of course.

First, we let $v_{k} := u - u_{k}$ and we observe that $v_{k}$ satisfies the Dirichlet boundary value problem
\begin{equation}
\left\{ \begin{array}{ll}
\L v_{k}=f-f(0),\quad \textit{in $\H_{k}$}\\
v_{k}=0,\quad \textit{in $\partial \H_{k}$}
\end{array} \right.
\end{equation}\label{dirichlet2}
From Lemma \ref{lemma-max} it follows that
\begin{align}\label{maxconse1}
\|v_{k}\|_{\infty}\le 4 \varrho^{2k}\|f-f(0)\|_{\infty}\le 4 \varrho^{2k} \omega_f(\varrho^k).
\end{align}
Moreover, since $\L (u_k-u_{k+1})=0$  in $\H_{k+1}$, we apply Proposition \ref{corollary} and Lemma \ref{lemma-max}, and we find
\begin{align}\nonumber
\Vert X_{I_i}(u_{k}-u_{k+1}) \Vert_{L^{\infty}(\H_{k+2})} &\leq C\varrho^{-k-2}\sup_{\H_{k+1}}|u_{k}-u_{k+1}|\nonumber \\ \nonumber
&\leq C \varrho^{-k}\Big(\sup_{\H_{k+1}}|v_{k}|+\sup_{\H_{k+1}}|v_{k+1}|\Big)\\
&\leq C \varrho^{-k}\varrho^{2k}\omega_{f}(\varrho^{k})=C\varrho^{k}\omega_{f}(\varrho^{k}), \label{terza}
\end{align}
for any $i= 1, \dots, m$. In the same way, we obtain
\begin{align}\nonumber
\Vert X_{I_i}X_{I_j}(u_{k}-u_{k+1}) \Vert_{L^{\infty}(\H_{k+2})} &\leq C\varrho^{-2k-4}\sup_{\H_{k+1}}|u_{k}-u_{k+1}| \\ 
&\leq C \varrho^{-2k}\varrho^{2k}\omega_{f}(\varrho^{k})=C\omega_{f}(\varrho^{k}) \label{quar}
\end{align}
for $i,j= 1, \dots, m$, and
\begin{align}\nonumber
\Vert X_0 (u_{k}-u_{k+1}) \Vert_{L^{\infty}(\H_{k+2})} &\leq C\varrho^{-2k-4}\sup_{\H_{k+1}}|u_{k}-u_{k+1}| \\ 
&\leq C \varrho^{-2k}\varrho^{2k}\omega_{f}(\varrho^{k})=C\omega_{f}(\varrho^{k}). \label{quarta}
\end{align}
Let $k\ge 1$ such that $\varrho^{k+4} \le \Vert z\Vert_\L \le \varrho^{k+3}$ , then we have:
\begin{equation}\label{stima_I_2}
\sum_{l=k}^{\infty} |\partial^2 u_l (0) - \partial^2 u_{l+1} (0)| \le C\sum_{l=k}^{\infty} \omega_f(\varrho^l)\le C \int_0^{\Vert z\Vert_\L} \frac{\omega_f (r)}{r} dr.
\end{equation}

We next identify the sum of the series $\sum_{l=k}^{\infty} \left(\partial^2 u_l (0) - \partial^2 u_{l+1} (0)\right)$ as 
\begin{equation}\label{sum-series}
 \sum_{l=k}^{\infty} \left(\partial^2 u_l (0) - \partial^2 u_{l+1} (0)\right) = \partial^2 u_k (0) - \partial^2 u (0).
\end{equation}
To this end, we first consider the derivative $ X_{I_i} X_{I_j}u_k$ and we prove that
\begin{eqnarray} \label{conv-Taylor}
\lim_{k \to +\infty} X_{I_i} X_{I_j} u_k(0)=X_{I_i} X_{I_j}T_0^2u(0), 
\end{eqnarray}
where $T_0^2u(\z)$ is the second-order Taylor polynomial of $u$ around the origin, computed at some point $\z\in \H_k$:
\begin{eqnarray*}
T_0^2u (\z)= u(0) + \sum_{deg(X_{I_i} \leq 2} X_{I_i}u(0) e_{i}
+ \frac{1}{2}\sum_{i,j=1}^m X_{I_i} X_{I_j}u(0) e_{i} e_{j} .
\end{eqnarray*}

Let us explicitly note that from \eqref{JacobianE} it follows that 
$${X_i h_j} _{\{h=0\}} = \delta_{ij},$$
where $\delta_{ij}$ is the Kronecker Delta function. 
Thus, by applying Theorem \ref{taylor} to $u\in C^{2}_\L(\H_{1}(0))$, 
we obtain from \eqref{conv-Taylor} that
\begin{eqnarray} \label{conv-u}
\lim_{k \to +\infty} X_{I_i} X_{I_j} u_k(0)=X_{I_i} X_{I_j} u(0).
\end{eqnarray}

We compute $\L T_0^2 u$ in $\zeta=(\xi,\tau)$ as 
\begin{equation*}
\begin{split}
\L T_0^2 u(\zeta)& = \! \! \sum_{i,j=1}^m X_{I_i} X_{I_j}u(0) - X_0 u(0) + \sum_{deg(X_{I_j} = 2} X_{I_j} u(0)= f(0),
\end{split}
\end{equation*}
where the last equality directly follows from the fact that $u$ is a solution. 


We now apply Proposition \ref{corollary} to $T_0^2 u-u_k$ for $R=\varrho^k$ and infer
\begin{equation}\label{tay-uk}
 |\partial^2_{x_i x_j}(u_k - T_0^2 u)(0)| \leq C \varrho^{-2k}\sup_{\H_k} \vert u_k - T_0^2 u \vert .
\end{equation}

Moreover, since $T_0^2u$ is the second-order Taylor polynomial of $u$, we have $u(\z) = T_0^2u(\z) + o(\Vert \z \Vert_\L^2)$. It follows that 
\begin{eqnarray}\label{esti-with-tay}
\sup_{\z \in \H_k}|u - T_0^2u| = o(\varrho^{2k})
\end{eqnarray}
Thus, from estimates \eqref{esti-with-tay} and \eqref{maxconse1}, we obtain
\begin{eqnarray}\label{sup-est}
\sup_{\H_k}|u_k - T_0^2u| \leq \sup_{\H_k}|v_k| +\sup_{\H_k}|u - T_0^2u| \leq 4\omega_{f}(\varrho^k)\varrho^{2k} + o(\varrho^{2k}) \leq o(\varrho^{2k}).
\end{eqnarray}
Estimates \eqref{tay-uk} and \eqref{sup-est} finally yield
\begin{equation*}
|X_{I_i} X_{I_j}(u_k - T_0^2u)(0)| \leq C \varrho^{-2k} \sup_{\H_k}|u_k - T_0^2u| \leq C \varrho^{-2k}o(\varrho^{2k})  \leq o(1),
\end{equation*}
where, as usual, the indexes $i$ and $j$ range from $1$ to $m$. Thus, for any $i,j=1,\ldots,m$ we have showed that \eqref{conv-Taylor} holds true. Repeating the same argument for the vector field $X_0$, and using again Theorem \ref{taylor}, we obtain:
\begin{eqnarray*}
\lim_{k \to +\infty}X_0 u_k(0)=X_0 T_0^2u(0)=X_0 u(0).
\end{eqnarray*}
In conclusion, using \eqref{stima_I_2}, we obtain:
\begin{equation}\label{stima_I_2-b}
I_2 \le \sum_{l=k}^{\infty} |\partial^2 u_l (0) - \partial^2 u_{l+1} (0)| \le C \int_0^{\Vert z\Vert_\L} \frac{\omega_f (r)}{r} dr,
\end{equation}
for $k\ge 1$ such that $\varrho^{k+4} \le \Vert z\Vert_\L \le \varrho^{k+3}$.
Similarly, we can estimate $I_3$ through the solution of $\L v=f(z)$ in $\H_j(z)$ and $v=u$ on $\partial \H_j(z)$ and obtain
\begin{equation}\label{stima_I_3}
I_3 \le \sum_{l=k}^{\infty} |\partial^2 u_l (z) - \partial^2 u_{l+1} (z)| \le C \int_0^{\Vert z\Vert_\L} \frac{\omega_f (r)}{r} dr.
\end{equation}

Finally, let us estimate $I_1$. 
Since $h_k=u_k-u_{k+1} \in C^\infty(\H_{k+2})$,  we can apply Proposition \ref{mean-value-lem} to the functions $X_{I_i} X_{I_j} h_k$ and $X_0 h_k$:
\en{
|X_{I_i} X_{I_j} h_k(z)-X_{I_i} X_{I_j} h_k(0)| \leq \frac{C}{\varrho^{k}}\Vert z \Vert_\L \Vert X_{I_i} X_{I_j}{ h_k} \Vert_{L^{\infty}(\H_{k+1})}
}
and 
\en{
|X_0 h_k(z)-X_0 h_k(0)| \leq \frac{C}{\varrho^{k}}\Vert z \Vert_\L \Vert X_0 h_k \Vert_{L^{\infty}(\H_{k+1})},
}
for $i,j=1,\ldots,m$. We can now apply once again \eqref{quar} to obtain
\en{
|X_{I_i}X_{I_j} h_k(z)-X_{I_i}X_{I_j} h_k(0)| \leq \frac{C}{\varrho^{k}}\Vert z \Vert_\L \Vert X_{I_i}X_{I_j} h_k \Vert_{L^{\infty}(\H_{k+1})}
\leq C\Vert z \Vert_\L\varrho^{-k}\omega_{f}(\varrho^{k}).
}
In addition, thanks to \eqref{quarta}, we infer
\en{
|X_0 h_k(z)-X_0 h_k(0)| \leq \frac{C}{\varrho^{k}}\Vert z \Vert_\L \Vert X_0 h_k \Vert_{L^{\infty}(\H_{k+1})}
\leq C\Vert z \Vert_\L\varrho^{-k}\omega_{f}(\varrho^{k}) .
}
Hence, since $u_k(z)-u_k(0)=u_0(z)-u_0(0)+\sum_{j=0}^{k-1}\left(h_j(0)-h_j(z)\right)$, we have
\begin{align*}
I_1 &\le | \partial^2 u_0 (z)-\partial^2u_0 (0)|+\sum_{j=0}^{k-1} |\partial^2 h_j (z)-\partial^2 h_j(0)|\\
&\le C \Vert z \Vert_\L \big( \|u_0\|_{L^\infty(\H_0)}  +C\sum_{j=0}^{k-1} \varrho^{-j} \omega_f(\varrho^j)\big)\\
&\le C\Vert z \Vert_\L \big( \|u\|_{L^\infty(\H_1(0))} + \|f\|_{L^\infty(\H_1(0))}+ C \int_{\Vert z \Vert_\L} ^1 \frac{\omega_f(r)}{r^2}\big).
\end{align*}
Combining the above estimate with \eqref{stima_I_2-b} and \eqref{stima_I_3}, we complete the proof of \emph{(ii)}.

\medskip

 We consider $u_1$ solution to the following Dirichlet problem
\begin{eqnarray*}
\left\{ \begin{array}{ll}
\L u_{1}=f(0),\quad &\textit{in $\H_{1/2}(0)$}\\
u_{1}=u,\quad &\textit{in $\partial \H_{1/2}(0)$}
\end{array} \right.
\end{eqnarray*}
Then, we have
\begin{eqnarray}\label{par-u-0}
|\partial^2 u(0)| \leq |\partial^2 u(0) -\partial^2 u_1(0)| + |\partial^2 u_1(0)|.
\end{eqnarray}
Thanks to \eqref{stima_I_2-b}, we can estimate the first term in \eqref{par-u-0} as
\begin{eqnarray}
|\partial^2 u(0) -\partial^2 u_1(0)|\le C \int_0^{1} \frac{\omega_f (r)}{r} dr.
\end{eqnarray}
To estimate the second term in \eqref{par-u-0}, we consider the function $v(z):=u_1(z)\eta_{1/2}(z)$, where $\eta_{1/2} $ is the cut-off function introduced in \eqref{cut-off} with $R=\frac{1}{2}$. Reasoning as in the proof of Proposition \ref{lem-apriori}, we obtain
\begin{equation*}
\begin{split}
u(z) = v(z) &= \sum_{j= 1}^{m}\int_{\H_{\frac{1}{2}}(0)}\Big(\Gamma(z,\cdot) X_j ( X_j\eta_{1/2}) u_1\Big)(\z)d \z\\ 
&\quad- \int_{\H_{\frac{1}{2}}(0)}\Big(\Gamma(z,\cdot) X_0(\eta_{1/2})u_1\Big)(\z)d\z \\
&\quad - \int_{\H_{\frac{1}{2}}(0)}\Big(\Gamma(z,\cdot)\eta_{1/2}\L(u_1)\Big)(\z)d\z\\
&+ 2 \sum_{j= 1}^{m} \int_{\H_{\frac{1}{2}}(0)}\Big( X_j^\zeta \Gamma(z,\cdot)X_j \eta_{1/2}  u_1\Big)(\z)d\z,
\end{split}
\end{equation*}
where $z \in \H_{\frac{1}{4}}(0)$. Thanks to Lemma \ref{lemma-max}, we estimate
\begin{eqnarray*}
\sup_{\H_{\frac{1}{2}}(0)}{|u_1|} \leq \sup_{\H_{\frac{1}{2}}(0)}{|u|} +4|f(0)|. 
\end{eqnarray*}
As the derivatives of $\eta_{1/2}$ vanish in $\H_{3/8}(0)$, for any $i,j=1,\ldots,m$, we obtain
\begin{equation}\label{overlineI}
\begin{split}
|X_{I_i}X_{I_j}u_1(z)|
&\leq \sum_{j= 1}^{m}\int_{\H_{\frac{1}{2}}(0) \setminus \H_{\frac{3}{8}}(0)}\big\vert \Big(X_{I_i}X_{I_j} \Gamma(z,\cdot) (X_{j}X_{j} \eta_{1/2}) u_1\Big)(\z)\big\vert d\z \\ 
&\quad + \int_{\H_{\frac{1}{2}}(0) \setminus \H_{\frac{3}{8}}(0)}\big\vert\Big(X_{I_i}X_{I_j}\Gamma(z,\cdot) Y(\eta_{1/2}) u_1\Big)(\z)\big\vert d\z \\  
&\quad +2\sum_{j= 1}^{m} \int_{\H_{\frac{1}{2}}(0) \setminus \H_{\frac{3}{8}}(0)}\big\vert\Big(\langle X_{I_i}X_{I_j} X_j\zeta \Gamma(z,\cdot), X_j\eta_{1/2}\rangle u_1\Big) (\z)\big\vert d\z \\
&\quad +\bigg\vert f(0)\Big(X_{I_i}X_{I_j}\int_{\H_{\frac{1}{2}}(0)} [ \Gamma(z,\cdot) \eta_{1/2} ](\z) d\z \Big) \bigg\vert\\
&\quad =: \overline{I}_{1}(z) + \overline{I}_{2}(z)+\overline{I}_3(z)+\overline{I}_4(z).
\end{split}
\end{equation}
Moreover, as the derivatives of $\eta_{1/2}$ are bounded, we estimate the first and second integral in \eqref{overlineI} as
\begin{eqnarray*}
\overline{I}_{1}(z) \leq C \big[ \sup_{\H_{\frac{1}{2}}(0)}{|u|} +4|f(0)|\big],\\ 
\overline{I}_{2}(z) \leq C \big[ \sup_{\H_{\frac{1}{2}}(0)}{|u|} +4|f(0)|\big],\\
\overline{I}_{3}(z) \leq C \big[ \sup_{\H_{\frac{1}{2}}(0)}{|u|} +4|f(0)|\big].
\end{eqnarray*}
Finally, by taking advantage of \eqref{est-sec-gamma-const}, we obtain that $\overline{I}_4(z)$ is bounded by a constant $C$ that only depends on $B$, $\lambda$ and $\Lambda$.

By using the same argument we can estimate $|X_0 u_1(0)|$ and thus
\begin{eqnarray}\label{par-u-0-2}
|\partial^2 u_1(0)| \leq C \big[ \sup_{\H_{\frac{1}{2}}(0)}{|u|} +4|f(0)|\big].
\end{eqnarray}
Combining estimates \eqref{par-u-0} and \eqref{par-u-0-2}, we conclude the proof of Theorem \ref{th-2} with constant coefficients. 
}

\setcounter{equation}{0}\setcounter{theorem}{0}
\section{H\"older estimates for  operators with Dini continuous coefficients}
This Section is devoted to the proof of Theorem \ref{th-2}. We therefore consider a solution $u$ to the equation
\begin{equation*} 
\LL u=f,
\end{equation*}
where the operator $\LL$ does satisfy the hypotheses {\rm \bf[H.1]} and {\rm \bf[H.2]} and $f$ is assumed to be Dini continuous, and we proceed as in in the case with constant coefficients. Specifically, we denote $\H_{k}=\H_{\varrho^{k}}(0)$, $\varrho=\frac{1}{2}$ and we consider the following sequence of Dirichlet problems:
\begin{equation} \label{dirichlet3}
\left\{ \begin{array}{ll}
\sum\limits_{i,j=1}^m a_{ij}(0,0) X_{I_i}X_{I_j} u_{k}+X_0u_{k}=f(0), \quad \text{in $\H_{k}$}\\
u_{k}=u,\quad \text{on $\partial \H_{k}$}.
\end{array} \right.
\end{equation}
Note that the bounds given in Propositions \ref{lem-apriori}, \ref{corollary} and \ref{mean-value-lem} only depend on the constants $\lambda, \Lambda$ in {\rm \bf[H.2]} and on the matrix $B$. Keeping in mind this fact, the proof of Theorem \ref{th-2} is given by the same argument used in the case with constant coefficients.

\begin{proof}[Proof of Theorem \ref{th-2}]
Consider, for every $k \in \N$, the auxiliary function $v_{k} := u - u_{k}$, and note that it is a solution to the boundary value problem
\begin{equation} \label{dirichlet4}
\left\{ \begin{array}{lll}
\sum\limits_{i,j=1}^m a_{ij}(0,0) X_{I_i}X_{I_j} v_k +X_0 v_k \\
\quad =f-f(0)+\sum\limits_{i,j=1}^m (a_{ij}(0)-a_{ij}(x,t))X_{I_i}X_{I_j}  u, \quad \text{in $\H_{k}$}\\
v_k =0,\quad \textit{in $\partial \H_{k}$}
\end{array} \right.
\end{equation}
In order to simplify the notation, we let
\begin{equation} \label{eta-def}
    \eta := \max_{i,j=1,\dots,m}\|X_{I_i}X_{I_j} u\|_{L^\infty(\H_{1})}.
\end{equation}
From Lemma \ref{lemma-max} it follows that
\begin{align*}
\|v_{k}\|_{L^\infty(\H_k)} \le C\varrho^{2k} [\omega_f(\varrho^k)+\omega_a(\varrho^k) \eta].
\end{align*}
Hence
\begin{align*}
\|u_{k}-u_{k+1}\|_{L^\infty(\H_{k+1})} \le C\varrho^{2k} [\omega_f(\varrho^k)+\omega_a(\varrho^k) \eta].
\end{align*}

As already observed, we can apply Corollary \ref{corollary} and obtain estimates for the second order derivatives of $v_k$. In fact, for any $i,j=1,\ldots,m$, we have

\begin{align}\nonumber
\Vert X_{I_i}X_{I_j}(u_{k}-u_{k+1}) \Vert_{L^{\infty}(\H_{k+2})} &\leq C(\varrho^{k})^{-2}\sup_{\H_{k+1}}|u_{k}-u_{k+1}| \\ 
&\leq C \varrho^{-2k}\varrho^{2k}[\omega_{f}(\varrho^{k})+\omega_a(\varrho^k) \eta]=C[\omega_{f}(\varrho^{k}) +\omega_a(\varrho^k) \eta]
\end{align}
and
\begin{align}\nonumber
\Vert X_0(u_{k}-u_{k+1}) \Vert_{L^{\infty}(\H_{k+2})} &\leq C(\varrho^{k})^{-2}\sup_{\H_{k+1}}|u_{k}-u_{k+1}| \\ 
&\leq C \varrho^{-2k}\varrho^{2k}[\omega_{f}(\varrho^{k})+\omega_a(\varrho^k) \eta]=C[\omega_{f}(\varrho^{k}) +\omega_a(\varrho^k) \eta]
\end{align}

To estimate the second order derivatives of the function $u$, we apply Theorem \ref{taylor} and proceed as in the in the case with constant coefficients. Since there are no significant differences, we omit the details here.
\end{proof}

\end{document}